\documentclass[11pt]{article}
\usepackage{amsmath}
\usepackage{mathrsfs}
\usepackage{amsfonts}

\usepackage{amsfonts, amsmath, amssymb}
\usepackage{amssymb,amsfonts,amsmath,color,
latexsym, epsfig,cite, psfrag,eepic,colordvi}
\usepackage{amscd,graphics}
\usepackage{lineno}
\usepackage{comment}
\usepackage{verbatim}

\newcommand{\qed}{\hfill $\Box $}
\newcommand{\pf}{\noindent {\bf Proof.} }

\newtheorem{theorem}{Theorem}[section]
\newtheorem{lemma}[theorem]{Lemma}
\newtheorem{coro}[theorem]{Corollary}

\newtheorem{conjecture}[theorem]{Conjecture}

\begin{document}

\title{On stability of rainbow matchings}

\author{Hongliang Lu\footnote{Partially supported by National Natural
Science Foundation of China under grant No. 11871391}\\
School of Mathematics and Statistics\\
Xi'an Jiaotong University\\
Xi'an, Shaanxi 710049, China\\
\medskip \\
Yan Wang\footnote{Partially supported by National Key R\&D Program of China under grant No. 2022YFA1006400, National Natural Science Foundation of China under grant No. 12201400 and Explore X project of Shanghai Jiao Tong University}\\
School of Mathematical Sciences, CMA-Shanghai\\
Shanghai Jiao Tong University\\
Shanghai 200240, China\\
\medskip \\
Xingxing Yu\footnote{Partially supported by NSF grant DMS-1954134}\\
School of Mathematics\\
Georgia Institute of Technology\\
Atlanta, GA 30332, USA}

\date{}

\maketitle

\begin{abstract}


We show that for any integer $k\ge 1$ there exists an integer $t_0(k)$
such that for integers $t, k_1, \ldots, k_{t+1}, n$ with $t>t_0(k)$,
$\max\{k_1, \ldots, k_{t+1}\}\le k$,  and $n > 2k(t+1)$, the following
holds:
If $F_i \subseteq {[n]\choose k_i}$ and $|F_i|> {n\choose k_i}-{n-t\choose k_i} - {n-t-k \choose k_i-1} + 1$
for all $i \in [t+1]$, then either $\{F_1,\ldots, F_{t+1}\}$ admits a
rainbow matching of size $t+1$ or there exists $W\in {[n]\choose t}$
such that $W$ is a vertex cover of $F_i$ for all $i\in [t+1]$.
This may be viewed as a rainbow non-uniform extension of the classical Hilton-Milner theorem.
We also show that the same holds for every $t$ and $n > 2k^3t$,
generalizing a recent stability result of Frankl and Kupavskii on matchings to
rainbow matchings.
\end{abstract}

\section{Introduction}

For any positive integer $n$, let $[n]$ denote the  set $\{1,\ldots,n\}$.
For any nonnegative integer $k$ and any set $S$, let ${S\choose
  k}=\{e\subseteq S: |e|=k\}$.
Let $k \ge 2$ be an integer.
A \emph{$k$-uniform hypergraph} or \textit{$k$-graph}  is a pair
$H=(V,E)$, where $V=V(H)$ is the set of {\it vertices} of $H$ and
 $E=E(H)\subseteq {V\choose k}$ is the set of {\it edges} of $H$, and
 we often identify $H$ with $E(H)$ and use $e(H)$ or $|H|$ to denote the number of edges in $H$.
For any $S \subseteq V(H)$, let $H[S]$ denote the
subgraph of $H$ with $V(H[S])=S$ and $E(H[S])=\{e\in E(H): e \subseteq
S\}$, and let $H-S:= H[V(H)\setminus S]$.
Let $\bar{H}$ be the $k$-uniform complement of $H$, i.e.,
$V(\bar{H})=V(H)$ and $E(\bar{H})={V(H)\choose k}\setminus E(H)$.

A \emph{matching} in a hypergraph $H$ is a subset of $E(H)$ consisting
of pairwise disjoint edges. The maximum size of a matching in a hypergraph $H$
is denoted by $\nu(H)$. A {\it vertex cover} of $H$ is a set of vertices of $H$ that meets all edges of $H$, and we use $\tau(H)$ to denote the minimum size of a vertex cover of $H$.
Clearly, $\nu(H)\le \tau(H)$.

A classical problem in extremal set
theory is to determine $\max e(H)$ with $\nu(H)$ fixed.
Erd\H{o}s
\cite{Erdos65} in 1965 made the following conjecture:
For positive integers $k,n,t$ with $n \ge kt$, every $k$-graph $H$ on $n$ vertices with $\nu(H) <
t$ satisfies $e(H)\leq \max \left\{ {n\choose k}-{n-t+1\choose k}, {kt-1\choose k} \right\}.$
This bound is tight for the complete $k$-graph on $kt-1$ vertices and for the $k$-graph on $n$ vertices in which every edge intersects a fixed set of $t-1$ vertices.
There have been recent activities on this conjecture, see \cite{AHS12,AFH12,FLM,Fr13,Fr17,HLS,LM}.
In particular,  Frankl \cite{Fr13}
proved that if  $n\geq (2t-1)k-(t-1)$ and $\nu(H)<t$ then $e(H)\le
{n\choose k}-{n-t+1\choose k}$, with further improvement  by Frankl
and Kupavskii \cite{FK18}.

The work in this paper was motivated by a recent result of Frankl and Kupavskii \cite{FK19} on a stability version  of the Erd\H{o}s Matching Conjecture.
Let $n,k,s$ be positive integers. For $i\in [k]$ define
$${\cal A}_i^{(k)}(n,s):=\left\{A\in {[n]\choose k} : \left| A\cap [(s+1)i-1] \right| \ge i\right \}, $$
and let ${\cal H}^{(k)}(n,s)$ denote the family
$$\left\{e\in {[n]\choose k}: e\cap [s-1]\ne \emptyset\right\} \cup \left\{ [s+k]\setminus [s] \right\} \cup \left\{e\in {[n]\choose k}: e\cap [s]=\{s\}, e\cap ([s+k]\setminus [s])\neq\emptyset\right\}.$$
Note that $e({\cal A}_i^{(k)}(n,s))\le {n\choose k}
-\sum_{j=0}^{i-1}{(s+1)i-1 \choose j}$, $\nu({\cal H}^{(k)}(n,s))=s$, and   $e({\cal H}^{(k)}(n,s))={n\choose k}-{n-s\choose k}-{n-s-k\choose k-1}+1$.

Frankl and Kupavskii \cite{FK19} proposed the following conjecture.
 \begin{conjecture}\label{FK-stablity}
Let $n,k,t$ be three positive integers such that $n>kt$.  Let $F$ be a
$k$-graph of order $n$. If $\nu(F)\leq t$, then $\tau(F)\le t$ or
\[
e(F)\leq \max\{ e({\cal A}_2^{(k)}(n,t)), \dots, e({\cal A}_k^{(k)}(n,t)), e({\cal H}^{(k)}(n,t)) \}.
\]
\end{conjecture}

The case $t=1$ is the Hilton-Milner theorem in \cite{HM67}.
In  \cite{FK19}, Frankl and Kupavskii confirms Conjecture
\ref{FK-stablity} for $n\geq (2+o(1))kt$. We are interested in
extending results like this to rainbow matchings.

Let $\mathcal{F} = \{F_1,\ldots, F_q\}$ be a family
of hypergraphs. A set of pairwise disjoint edges, one from
each $F_i$, is called a \emph{rainbow matching} for $\mathcal{F}$, and we say that ${\cal F}$ or
$\{ F_1,\ldots, F_q\}$ {\it admits} a rainbow matching. Aharoni and
Howard \cite{AH} (also see Huang, Loh and Sudakov \cite{HLS})
conjectured the following rainbow version of the Erd\H{o}s Matching
Conjecture: If $F_i\subseteq {[n]\choose k}$ and $|F_i|> \max\left\{{n\choose k}-{n-q+1\choose k},{kq-1\choose k}\right\}$ for $i\in [q]$ then
$\{ F_1,\ldots, F_q\}$ admits a rainbow matching. This conjecture was confirmed for $n>3k^2q$ by Huang, Loh and Sudakov \cite{HLS}. The bound $3k^2q$ on $n$ was improved to $12kq \log(e^2q)$
by Frankl and Kupavskii \cite{FK20}.  Lu, Wang and Yu \cite{LWY20} further improved this bound to $n > 2kq$ for $q$ sufficiently large.
Recently, Keevash, Lifshitz, Long and Minzer \cite{KLLM2} proved a more general result with $n = \Omega(kq)$ using sharp threshold techniques developed in \cite{KLLM1}.

Our main result in this paper is the following, which extends the result of Frankl and
Kupavskii \cite{FK19} to rainbow matchings.

\begin{theorem}\label{main}
Let $k\ge 1$ be any integer. There exists  $t_0=t_0(k)$ such that, for
positive integers $t, k_1,k_2,\ldots,k_{t+1}, n$ with $t>t_0$,
$\max\{k_1, \ldots, k_{t+1}\}\le k$, and $n > 2k(t+1)$, the following
holds: If
$F_i \subseteq {[n]\choose k_i}$ and
\[
|F_i|> {n\choose k_i}-{n-t\choose k_i} - {n-t-k \choose k_i-1} + 1
\]
for all $i \in [t+1]$, then

\begin{itemize}
  \item [(i)] $\{F_1, \ldots, F_{t+1}\}$ admits a rainbow matching of size
    $t+1$; or
 \item [(ii)]  there exists $W \in {[n]\choose t}$  such that $W$ is a
   vertex cover of $F_i$ for all $i \in [t+1]$.
   \end{itemize}
\end{theorem}



To prove Theorem~\ref{main}, we need to establish the following
result, which may be viewed as a stability version(as well as a
generalization)  of  the above mentioned result of Huang, Loh and
Sudakov \cite{HLS} as well as an earlier result of Bollob\'as, Daykin and Erd\H{o}s \cite{BDE76}.

\begin{theorem}
\label{thm:stability-non-uniform}
Let $n, t, k_1,k_2,\ldots,k_{t+1}$ be positive integers such that $k_1 \ge k_2 \ge \ldots \ge k_{t+1} \ge 2$ and $n \ge 2 k_1^2 k_2 t$.
Let $F_i \in {[n] \choose k_i}$, $i \in [t+1]$, such that, for $i\in [2]$,
$$|F_i| > {n \choose k_i} - {n-t \choose k_i} - {n-t-k_{3-i} \choose k_i-1} + 1,$$
and for $i\in [t+1]\setminus [2]$,
$$|F_i| > {n \choose k_i} - {n-t \choose k_i} - {n-t-k_{3} \choose k_i-1} + 1.$$
Then one of the following holds:
\vspace{-4mm}
\begin{itemize}
	 	\setlength{\itemsep}{0pt}
		\setlength{\parsep}{0pt}
		\setlength{\parskip}{0pt}
\item[(i)] $\{F_1, \ldots, F_{t+1}\}$ admits a rainbow matching of size $t+1$; or
\item[(ii)] there exists $W \in { [n]\choose t}$ such that $W$ is a vertex cover of $F_i$ for all $i \in [t+1]$.
\end{itemize}
\end{theorem}


In Section 2, we first reduce the problem for finding rainbow matchings for non-uniform hypergraphs to one for uniform hypergraphs (through an operation called ``expansion''), and then reduce the problem to
a matching problem for ${\cal F}^{t+1}(k,n)$, a special class of
uniform hypergraphs. We will see that the extremal hypergraphs of Theorem \ref{main} naturally corresponds to a special class of
$(k+1)$-graphs,  denoted by ${\cal F}_{t+1}(n,k;t)$ and defined later.
This matchings problem for ${\cal F}^{t+1}(k,n)$ will be further reduced to a near perfect matching problem.

In Section 3, we prove Theorem \ref{thm:stability-non-uniform}  by an inductive argument.
We prove Theorem~\ref{main} in Section 4 when ${\cal F}^{t+1}(k,n)$ is close to ${\cal F}_{t+1}(k,n;t)$, in the sense that most edges of ${\cal F}_{t+1}(k,n;t)$ are also edges of ${\cal F}^{t+1}(k,n)$.
In Section 5, we deal with the case when ${\cal F}^{t+1}(k,n)$ is not close to ${\cal  F}_{t+1}(k,n;t)$. We conclude with a rainbow version of Conjecture~\ref{FK-stablity} in Section 6.


\section{Reductions}
The goal of this section is to reduce the problem for finding  rainbow matchings for
non-uniform hypergraphs to  a near perfect matching problem for
uniform hypergraphs.
We need to use shadows of hypergraphs.
Let $k, m$ be positive integers.
The {\it $k$-cascade representation} of $m$ is
$$m = {a_k \choose k} + {a_{k-1} \choose k-1} + \cdots + {a_s \choose s}$$
where $a_k > a_{k-1} > \dots > a_s \ge s \ge 1$ are integers.
Given a family ${\cal F}$ of sets, the {\it shadow} $\partial {\cal F}$ of
${\cal F}$ is defined as
$$\partial {\cal F} = \{E : E = F \setminus \{x\} \text{ for some } F \in {\cal F} \text{ and } x \in F \}.$$
For $i \ge 1$, we define $\partial^{i+1} {\cal F} = \partial (\partial^{i} {\cal F})$.
 The following result is known as the Kruskal-Katona theorem (see Theorem 1 in \cite{Fr84}\label{KK}).

\begin{theorem}[Kruskal and Katona]\label{KK}
Let $n,k$ be positive integers and ${\cal F} \subseteq {[n] \choose k}$.
If $|{\cal F}| = {a_k \choose k} + {a_{k-1} \choose k-1} + \cdots + {a_s \choose s}$ where $a_k > a_{k-1} > \dots > a_s \ge s \ge 1$ are integers.
Then
\[
|\partial {\cal F}| \ge {a_k \choose k-1} + {a_{k-1} \choose k-2} + \cdots + {a_s \choose s-1}.
\]
\end{theorem}

The first step of our reduction needs the following connection between the number of edges in an $r$-graph and its $k$-uniform ``expansion".

\begin{lemma}\label{Obser3}
Let $n,k,r$ be three positive integers with $n\ge k>r$. Let $F\subseteq {[n]\choose r}$ and let $G= \{e \cup f: e \in F, f \in {[n] \setminus e \choose k-r} \}.$
If $|F|>{n\choose r}-{n-t\choose r} - {n-t-k\choose r-1} + 1$, then $|G|>{n\choose k}-{n-t\choose k} - {n-t-k \choose k-1} + 1$.
\end{lemma}

\pf Suppose the conclusion is false. Then the number of edges in $\overline{G}$, the complement of $G$, satisfies $|\overline{G}| \ge {n-t\choose k} + {n-t-k \choose k-1} - 1$. Applying the
formula ${x\choose y}={x-1\choose y}+{x-1\choose y-1}$ to ${n-t-k\choose k-1}$ and repeat, we obtain
\begin{equation*}
\begin{split}
|\overline{G}|
&\geq {n-t\choose k} + {n-t-k-1 \choose k-1} + {n-t-k-1 \choose k-2} - 1 \\
&= {n-t\choose k} + {n-t-k-1 \choose k-1} + {n-t-k-2 \choose k-2} + {n-t-k-2 \choose k-3} - 1 \\
&= {n-t\choose k} + \left(\sum_{i=1}^{k-1}{n-t-k-i \choose k-i} +1\right) -1\\
&= {n-t\choose k} +\sum_{i=1}^{k-1}{n-t-k-i \choose k-i}
\end{split}
\end{equation*}
Thus, by Theorem \ref{KK}, we have
\begin{equation*}
\begin{split}
|\partial^{k-r} (\overline{G})|
&\ge {n-t\choose r} +\sum_{i=1}^{r}{n-t-k-i \choose r-i} -1\\
&=  {n-t\choose r} + {n-t-k \choose r-1} - 1
\end{split}
\end{equation*}
Since $\partial^{k-r} (\overline{G}) \cap F = \emptyset$ (by definition),  $|F| \le {n\choose r}-{n-t\choose r} - {n-t-k \choose r-1} + 1$, a contradiction.
\qed

\medskip

By Lemma \ref{Obser3}, Theorem \ref{main} follows from the following result.
\begin{theorem}\label{main_uniform}
Let $k \ge 3$ be a positive integer. There exists $t_0=t_0(k)$ such
that, for integers $t,n$ with  $t > t_0$ and $n > 2k(t+1)$, if $F_i \subseteq {[n]\choose k}$ and
\[
|F_i|> {n\choose k}-{n-t\choose k} - {n-t-k \choose k-1} + 1,
\]
for $i\in [t+1]$, then

\begin{itemize}
  \item [(i)] $\{F_1, \ldots, F_{t+1}\}$ admits a rainbow matching; or
   \item [(ii)] there exists $W \in {[n]\choose t}$ such that $W$ is a vertex cover of $F_i$ for all $i\in [t+1]$.
\end{itemize}
   \end{theorem}

To prove Theorem~\ref{main_uniform}, we convert
this rainbow matching problem on $k$-graphs  to a matching problem for a special class of $(k+1)$-graphs.
Let $F_1,\ldots, F_q$ be a family of subsets of  ${[n]\choose k}$ and $X:=\{x_1,\ldots,x_q\}$ be a set of $q$ vertices.
We use $\mathcal{F}^q(k,n)$ to denote the $(k+1)$-graph  with  vertex
set  $X\cup  [n]$ and edge set
\[
E(\mathcal{F}^q(k,n))=\bigcup_{i=1}^q \{\{x_i\}\cup e\ :\ e\in F_i\}.
\]
If $F_1=\cdots=F_q=H_k(s,n)$, where  $1\le s<n$ and $H_k(s,n)$ denotes the $k$-graph with vertex set $[n]$
 and edge set  ${[n]\choose k}\setminus {[n]\setminus [s]\choose k}$, then
we denote such $\mathcal{F}^q(k,n)$ by $\mathcal{F}_q(k,n;s)$.

\medskip

\textbf{Observation 1:}  $\{F_1,\ldots,F_q\}$ admits a rainbow matching if, and only if,  $\mathcal{F}^q(k,n)$ has a matching of size $q$.

\medskip

Hence, to prove Theorem~\ref{main_uniform}, we need to see when
$\mathcal{F}^{t+1}(k,n)$ has a matching of size $t+1$. We further
reduce this problem to a near perfect matching problem.

Write $n-kq=km+r$, where $0\leq r\leq k-1$. Let $F_1, \ldots, F_q\subseteq {[n]\choose k}$, and let $F_i={[n]\choose k}$ for $i=q+1, \ldots, q+m$.
  Let $Q=\{x_1,\ldots,x_{m+q}\}$ and let $\mathcal{H}^{q}(k,n)$ be
  the $(k+1)$-graph with vertex set $Q\cup [n]$ and edge set
\begin{align*}
E(\mathcal{H}^{q}(k,n))=\left( \bigcup_{i=1}^{q}\{\{x_i\}\cup e\ :\ e\in F_i\} \right) \bigcup \left( \bigcup_{i=q+1}^{q+m}\{\{x_i\}\cup e\ :\ e\in  {[n] \choose k} \} \right).
\end{align*}
When $F_1=\cdots=F_{q}=H_k(s,n)$,
we denote $\mathcal{H}^{q}(k,n)$ by $\mathcal{H}_{q}(k,n;s)$. Note that  $\nu(\mathcal{H}_{q}(k,n;q))=q+m=(n-r)/k$, i.e., ${\cal H}_q(k,n;q)$ has a matching covering
all but at most $r<k$ vertices (and such a matching is said to be \emph{near perfect}).

The following lemma provides an equivalent condition on matchings in $\mathcal{F}^{q}(k,n)$ and matchings in $\mathcal{H}_{q}(k,n;q)$.  See Lemma 2.1 in \cite{LWY20}.

\begin{lemma}\label{Rain-PM}
Let $n,k,q$ be positive integers and let $F_1, \ldots, F_{q}\subseteq {[n]\choose k}$. Then
$\mathcal{F}^{q}(k,n)$ has a matching of size $q$ if, and only if,  $\mathcal{H}_q(k,n;q)$ has a matching of size $m+q$, where $m+q=\lfloor n/k\rfloor$.
\end{lemma}



\medskip

For  the proof of Theorem \ref{main_uniform}, we need to tell how far
a hypergraph is close to another hypergraph.
Given two $k$-hypergraphs $H_1, H_2$ with $V(H_1)=V(H_2)$, let $c(H_1,H_2)$ be the minimum of $|E(H_1)\backslash E(H')|$
taken over all isomorphic copies $H'$ of $H_2$ with  $V(H') = V(H_2)$.
For a real number $\varepsilon > 0$,
we say that $H_2$ is \textit{$\varepsilon$-close} to $H_1$ if $V(H_1) = V(H_2)$ and $c(H_1,H_2)\leq \varepsilon|V(H_1)|^k$.
The following is striaghtforward to verify.

\noindent\textbf{Observation 2:}  (i) If $\mathcal{F}^q(k,n)$ is $\varepsilon$-close to $\mathcal{F}_q(k,n;s)$ then
$\mathcal{H}^q(k,n)$ is $\varepsilon$-close to $\mathcal{H}_q(k,n;s)$; (ii) for $n \le 2k^3\min\{q,s\}$, if $\mathcal{H}^q(k,n)$ is $\varepsilon$-close to $\mathcal{H}_q(k,n;s)$ then
$\mathcal{F}^q(k,n)$ is $4k^5 \varepsilon$-close to $\mathcal{F}_q(k,n;s)$.

Our proof of Theorem~\ref{main_uniform} will be divided into two parts, according to  whether or
not  $\mathcal{F}^{t+1}(k,n)$ is close  to $\mathcal{F}_{t+1}(n,k;t)$.
If $\mathcal{F}^{t+1}(k,n)$ is close to $\mathcal{F}_{t+1}(n,k;t)$, we will apply greedy argument to construct a matching of size $t$.
If $\mathcal{F}^{t+1}(k,n)$ is not close to $\mathcal{F}_{t+1}(n,k;t)$, then by Observation 2, $\mathcal{H}^{t+1}(k,n)$ is not close to $\mathcal{H}_{t+1}(n,k;t)$ (see proof of Theorem \ref{main_uniform}), and
we will  show that $\mathcal{H}^{t+1}(k,n;t)$ has a small
matching $M_1$ with nice absorbing properties and
$\mathcal{H}^{t+1}(k,n;t)-V(M_1)$ has a spanning subgraph
in which we can find a large matching $M_2$ that can be extended to a
near perfect matching by using $M_1$.

\section{Stability results on small matchings}

We begin with the result of Huang, Loh and Sudakov \cite{HLS}
mentioned previously, and state a corollary which generalizes it to
non-uniform hypergraphs.

\begin{theorem}[Huang, Loh and Sudakov]\label{HLS}
Let $n,k,t$ be three positive integers such that $n > 3k^2t$.
Let ${\cal F}=\{F_1,\ldots, F_t\}$ be a family of subsets of ${[n]\choose k}$. If
\[
|F_i|> {n\choose k}-{n-t+1\choose k}
\]
for all $1\leq i\leq t$, then ${\cal F}$ admits a rainbow matching.
\end{theorem}


\begin{coro}\label{thm:HLS-non-uniform}
Let $t,n,k_1,\ldots,k_t,k$ be positive integers such that $k = \max\{k_1,\ldots,k_t\}$ and $n > 3k^2t$.
Let ${\cal F}=\{F_1,\ldots, F_t\}$ be a family where $F_i \in {[n]\choose k_i}$ for $i \in [t]$. If
\[
|F_i|> {n\choose k_i}-{n-t+1\choose k_i}
\]
for all $1\leq i\leq t$, then ${\cal F}$ admits a rainbow matching.
\end{coro}

\pf
For $i \in [t]$, we define the expansion of $F_i$ as
$$G_i = \{e \cup f: e \in F_i, f \in {[n] \setminus V(e) \choose k-k_i} \}.$$
Note that each $G_i$ is a $k$-graph.

We claim that  $|\bar{G_i}|< {n-t+1\choose k}$ for all $i\in [t]$. For, otherwise, there exists some $i\in [t]$ such that $|\bar{G_i}|\ge {n-t+1\choose k}$.
Then by Theorem \ref{KK}, $|\partial^{k-k_i} (\bar{G_i})| \ge {n-t+1\choose k_i}$.
Since $\partial^{k-k_i} (\bar{G_i}) \cap F_i = \emptyset$ (by definition), $|F_i| \le {n\choose k_i}-{n-t+1\choose k_i}$, a contradiction.

Therefore, $|G_i|> {n\choose k}-{n-t+1\choose k}$. Hence, by Theorem \ref{HLS}, $\{G_1,\ldots, G_t\}$ admits a rainbow matching, which implies that  ${\cal F}$ admits a rainbow matching.
\qed


%
%
\medskip

We also need a classical theorem of M\"{o}rs, which is Theorem 6 in \cite{Mo85}.

\begin{theorem}[M\"{o}rs, \cite{Mo85}]
\label{thm:mors}
Let $n,k,l$ be three positive integers such that $n \ge k+l$.
Let $\mathcal{A} \in {[n] \choose k}$ and $\mathcal{B} \in {[n] \choose l}$.
Suppose $\cap_{C \in \mathcal{A} \cup \mathcal{B}} C = \emptyset$ and
$A \cap B \ne \emptyset$ for all $A \in \mathcal{A}$ and $B \in \mathcal{B}$.
Then either
$$|\mathcal{A}| \le {n-1 \choose k-1} - {n-1-l \choose k-1} + 1$$
or
$$|\mathcal{B}| \le {n-1 \choose l-1} - {n-1-k \choose l-1} + 1.$$
\end{theorem}


\noindent {\bf Proof of Theorem \ref{thm:stability-non-uniform}.}
We apply induction on $t$. For the base case, suppose $t=1$. Then for
$i\in [2]$,
$$|F_i| >  {n\choose k_i}-{n-1\choose k_i}-{n-1-k_{3-i}\choose
  k_i-1}+1={n-1\choose k_i-1}-{n-1-k_{3-i}\choose k_i-1}+1;$$
so (i) of Theorem~\ref{thm:stability-non-uniform} follows from Theorem
\ref{thm:mors}.

Now suppose $t\ge 2$ and the conclusion holds with $t$ sets.
Moreover,
we may assume that $\{F_1, \ldots, F_{t+1}\}$ does not admit a rainbow matching, as otherwise (i) holds.

Since $k_1\ge k_i$ for $i\in [t+1]$, we have from the assumption of Theorem \ref{thm:stability-non-uniform} that $|F_i| \ge {n\choose k_i}-{n-t\choose k_i}-{n-t-k_1\choose
  k_i-1}+1$ for $i\in [t+1]$. Hence, since  ${n-(t-1)\choose k_i}={n-t\choose
  k_i}+{n-t\choose k_i-1}$, $|F_i|> {n \choose k_i} - {n-(t-1)
  \choose k_i}$ for $i \in [t]$. Thus by Corollary
\ref{thm:HLS-non-uniform},
$\{F_1,\ldots,F_t\}$ admits a rainbow matching, say  $M_1$.
Therefore, since $\{F_1, \ldots, F_{t+1}\}$ does not admit a rainbow matching,
every edge in $F_{t+1}$ must intersect $V(M_1)$.
Thus the maximum degree $\Delta(F_{t+1}) \ge |F_{t+1}| / |V(M_1)|$.
Let $v \in [n]$ such that $d_{F_{t+1}}(v) = \Delta(F_{t+1})$.
Then, since
$${n\choose k_{t+1}}={n-1\choose k_{t+1}-1}+{n-1\choose
  k_{t+1}} =\ldots =\sum_{i=1}^t{n-i\choose k_{t+1}-1} +{n-t\choose
  k_{t+1}}\ge t{n-t\choose k_{t+1}-1} +{n-t\choose k_{t+1}},$$
we have
$$d_{F_{t+1}}(v) \ge \frac{|F_{t+1}|}{|V(M_1)|} > \frac{(t-1){n-t \choose k_{t+1}-1}}{\sum_{i=1}^{t} k_i}.$$
Hence, since $n \ge 2 k_1^2 k_2 t$,
$$ d_{F_{t+1}}(v)  > \frac{(t-1){n-t \choose k_{t+1}-1}}{\sum_{i=1}^{t} k_i} > \left(\sum_{i=1}^{t} k_i \right) {n-2 \choose k_{t+1} - 2}.$$

If $\{F_1 - v, \ldots, F_t - v\}$ admits a rainbow matching,
say $M_2$, then $d_{F_{t+1}}(v)> |V(M_2)| {n-2 \choose k_{t+1} -
  2}$ as $|V(M_2)| = \sum_{i=1}^{t} k_i$. 
So there exists an edge $e \in F_{t+1}$ such that $v \in e$ and $e \cap V(M_2) =\emptyset$.
Now $M_2 \cup \{e\}$ is a desired rainbow matching for $\{F_1,
\ldots, F_t, F_{t+1}\}$, a contradiction.

So  $\{F_1 - v, \ldots, F_t - v\}$ does not admit any rainbow
matching.  Note that, for $i\in [2]$,
$$|F_i - v| \ge |F_i| - {n-1 \choose k_i - 1} > {n-1 \choose k_i} - {(n-1)-(t-1) \choose k_i} - {(n-1)-(t-1)-k_{3-i} \choose k_i-1} + 1,$$
and, for $i\in [t]\setminus [2]$,
$$|F_i - v| \ge |F_i| - {n-1 \choose k_i - 1} > {n-1 \choose k_i} - {(n-1)-(t-1) \choose k_i} - {(n-1)-(t-1)-k_{3} \choose k_i-1} + 1.$$
Thus, by inductive hypothesis (applied to $\{F_1 - v, \ldots, F_t - v\}$),
there exists $W' \in {[n] \setminus \{v\}\choose t-1}$ such that $W'$ is a vertex cover of $F_i-v$ for all $i \in [t]$.

Write $W = W' \cup \{v\}$. We may assume that there exists $f \in F_{t+1}$ such that $f \cap W = \emptyset$; for otherwise (ii) holds.
Note that the number of edges in $F_i$ (for each $i\in [t]$)  intersecting both $f$ and $W$
is
$$e_i(f, W) \le {n \choose k_i} - {n-k_{t+1} \choose k_i} -  {n-t \choose k_i} + {n-t-k_{t+1} \choose k_i}.$$
Hence,  for each $i \in [t]$, 
since every edge of $F_i$ intersects $W$,
the number of edges in $F_i$ disjoint
from $f$ is
\begin{equation*}
\begin{split}
|F_i - f|
&= |F_i| - e_i(f,W) \\
&>\left( {n \choose k_i} - {n-t \choose k_i} - {n-t-k_{3} \choose
    k_i-1} + 1 \right) -e(f,W) \\
&\ge {n-k_{t+1} \choose k_i} - {n-t-k_{3} \choose k_i-1} - {n-t-k_{t+1} \choose k_i} + 1 \\
&> {n-k_{t+1} \choose k_i} - {(n-k_{t+1})-(t-1) \choose k_i}
\end{split}
\end{equation*}
Hence by Corollary \ref{thm:HLS-non-uniform}, $\{F_1-f, F_2-f, \ldots, F_t-f\}$ admits a rainbow matching $M_3$ of size $t$.
Therefore, $M_3 \cup \{f\}$ is a rainbow matching of size $t+1$ which satisfies (i).
\qed

%
As a consequence of Theorem \ref{thm:stability-non-uniform}, we have the following conclusion.

\begin{coro}
\label{thm:exact2}
Let $t,  k_1,k_2,\ldots,k_{t+1}, n$ be positive integers, let  $k =
\max\{k_1, \ldots, k_{t+1}\}$, and let $F_i \subseteq {[n] \choose k_i}$ for $i \in [t+1]$.
Suppose  $n \ge 2 k^3t$ and $|F_i| > {n \choose k_i} - {n-t \choose k_i} - {n-t-k \choose k_i-1} + 1$ for all $i \in [t+1]$.
Then
\vspace{-4mm}
\begin{itemize}
	 	\setlength{\itemsep}{0pt}
		\setlength{\parsep}{0pt}
		\setlength{\parskip}{0pt}
\item[$(i)$] $ \{F_1,F_2,\cdots,F_{t+1}\}$ admits a rainbow matching
  of size $t+1$; or
\item[$(ii)$] there exists $W\in {[n]\choose t}$ such that $W$ is a vertex cover of $F_i$ for  all $i \in [t+1]$.
\end{itemize}
\end{coro}

\section{Extremal case}
%

In this section, we prove Theorem~\ref{main}  for the case when $\mathcal{F}^{t+1}({k},n)$
is $\varepsilon$-close to the extremal configuration $\mathcal{F}_{t+1}({k},n;t)$ for some $ \varepsilon\ll 1/k$.
We write $0<a\ll b$ to mean that there exists an increasing function
$f$ such that our result holds whenever $a\leq f(b) < b$.


Let $H$ be a $(k+1)$-graph and $v \in V(H)$. We define the neighborhood $N_H(v)$ of $v$ in $H$ to be the set $ \{ S \in {V(H) \choose k} \ :\ S \cup \{v\} \in E(H) \} $.
Let $H$ be a $(k+1)$-graph with the same vertex set as $\mathcal{F}_{t+1}(k,n;t)$.
Given real number $\alpha$ with $0<\alpha < 1$, a vertex $v$ in $H$ is called
\emph{$\alpha$-good} with respect to $\mathcal{F}_{t+1}(k,n;t)$ if
$$\left|N_{\mathcal{F}_{t+1}(k,n;t)}(v)\setminus N_H(v)\right|\le \alpha n^{k}.$$
Clearly, if $H$ is $\varepsilon$-close to $\mathcal{F}_{t+1}(k,n;t)$,
then at most $(k+1)(1+1/k)^{k+1}\varepsilon n/\alpha$ vertices of $H$ are not
$\alpha$-good with respect to $\mathcal{F}_{t+1}(k,n;t)$.

First we deal with the case when all vertices of $\mathcal{F}^{t+1}({k},n)$
are $\alpha$-good with respect to $\mathcal{F}_{t+1}(k,n;t)$.
Let $Q,V$ be two disjoint sets.  A $(k+1)$-graph $H$ with vertex $Q \cup V$
is called \emph{$(1,k)$-partite} with {\it partition classes} $Q,V$
if, for each edge $e\in E(H)$, $|e\cap Q|=1$ and $|e\cap
V|=k$. Clearly, $\mathcal{F}^{t+1}({k},n)$ and $\mathcal{F}_{t+1}({k},n;t)$ are $(1,k)$-partite graphs
with partition classes $X, [n]$.

\begin{lemma}\label{good-lem}
Let $\zeta, \alpha $ be real numbers and  $n,k,t$ be positive integers such that $0<\alpha \ll \zeta <1$, $\alpha\ll 1/k$, $n \ge 8 k^4$, $t \ge n/(2k^3)$ and $t+1<(1-\zeta)n/k$.
Let $H$ be a $(1,k)$-partite $(k+1)$-graph with the same partition
classes as $\mathcal{F}_{t+1}(k,n;t)$.
If every vertex of $H$ is $\alpha$-good  with respect to $\mathcal{F}_{t+1}(k,n;t)$, then $\nu(H)\geq t$ with equality only if $H$ is a subgraph of $\mathcal{F}_{t+1}(k,n;t)$.
\end{lemma}

\pf
Let $X := \{x_1, x_2, ..., x_{t+1}\}$, $W := [t]$, and $U := [n]
\setminus [t]$, such that $X,[n]$ are the partition classes of $H$ and
$\mathcal{F}_{t+1}(k,n;t)$.

If $H$ is a subgraph of $\mathcal{F}_{t+1}(k,n;t)$, then we define $e_0 = \emptyset$.
Otherwise, there exists an edge $e_0$ in $H$ such that $|e_0 \cap X| =
1$ and $e_0\cap W=e_0 \cap [t] = \emptyset$.
Without loss of generality, we may assume $e_0 \cap X = \{ x_{t+1} \}$.

Write $H' = H - V(e_0)$, i.e. deleting all vertices in $e_0$. Let $M$ be a maximum matching in $H'$ such that $|e\cap X| = |e\cap
W|=1$ for all $e\in M$. Thus $|M|\le |W|=t$.  Let $X'=X\setminus (V(M)
\cup V(e_0))$, $W'=W\setminus (V(M) \cup V(e_0))$, and $U'=U\setminus
(V(M) \cup V(e_0))$. Thus $|W'|=|W|-|M|$ and $|U'|\ge
|U|-|M|(k-1)-k$.

We claim that $|M|\geq n/(4k^3)$. For,
suppose $|M| <  n/(4k^3)$.
Consider any vertex $x\in X'$.
Since $x$ is $\alpha$-good with respect to $\mathcal{F}_{t+1}(k,n;t)$, we have
\[
\left| \left(W\times {U\choose k-1}\right)\setminus N_H(x)\right|\leq \alpha n^k.
\]
Since $t + 1 <(1-\zeta)n/k$, $|U'|\ge |U|-|M|(k-1)-k\ge
(n-t)-t(k-1)-k> \zeta n$. Hence,  since $t \ge n/2k^3$,
\[
\left|W'\times {U'\choose k-1}\right|\ge (|W|-|M|){|U|-|M|(k-1)-k \choose k-1}>\frac{n}{4k^3}{\zeta n\choose k-1}
>\frac{n}{4k^3} \frac{(\zeta n/2)^{k-1}}{(k-1)!},
\]
Therefore, $\left|W'\times {U'\choose k-1}\right| > \alpha n^k$, as $\alpha< \zeta^{k-1}(k^3 2^{k+1} (k-1)!)^{-1}$.
Hence,  there exists $f\in N_{H}(x)\cap \left(W'\times {U'\choose
    k-1}\right)$. Let $f'=\{x\}\cup f$, Then $f'\in E(H)$, $|f'\cap
X|=|f'\cap W|=1$, and $f'\cap V(M)=\emptyset$.  Now $M'=M\cup \{f'\}$ is a  matching of size $|M|+1$ in $H$, and
$|e\cap X|= |e\cap W|=1$ for all $e\in M'$. Thus, $M'$ contradicts the choice of $M$, completing the proof of the claim.

\medskip

Let $\{u_1,\ldots,u_{k+1}\}\subseteq V(H)\setminus V(M)$, where $u_1\in X'$, $u_{k+1}\in W'$ and $u_i\in U'$ for $i\in [k]\setminus \{1\}$.
Since $|M| \ge n/(4k^3) \ge 2k$, 
let $\{e_1,\ldots,e_{k}\}$ be an arbitrary $k$-subset of $M$, and let
$e_i := \{v_{i,1},v_{i,2},\ldots,v_{i,k+1}\}$ with $v_{i,1} \in X$, $v_{i,k+1}\in W$, and $v_{i,j}\in U$ for $i \in [k]$ and  $j \in [k] \setminus \{1\}$.
For $j \in [k+1]$, let $f_j := \{u_{j}, v_{1,j+1},
v_{2,j+2},\ldots,$ $v_{k,j+k}\}$ with addition in the subscripts
modulo $k+1$ (except we write $k+1$ instead of $0$). Note that $f_1, \ldots, f_{k+1}$ are pairwise disjoint.

If $f_j \in E(H)$ for all $j \in [k+1]$ then $M':= (M \cup
\{f_1,\ldots,f_{k+1}\})\setminus \{e_1,\ldots,e_{k}\}$ is a matching in $H$
such that  $|M'| = |M| + 1 > |M|$ and $|f\cap X|=|f\cap W|=1$ for all $f\in M'$, contradicting the choice of $M$.
Hence, $f_j\not\in E(H)$ for some $j \in [k+1]$.

Note that there are $\binom{|M|}{k}k!$  choices of ordered $k$-tuples $(e_1,\ldots, e_{k}) \in M^k$ and that any two different such choices correspond to different $f_j$.
Hence,
\begin{eqnarray*}
& & |\{e \in E(\mathcal{F}_{t+1}(k,n)) \setminus  E(H): |e\cap \{u_i: i\in [k+1]\}|=1\}|\\
&\geq & |M|(|M| - 1) \cdots (|M| - k +1) \\
&> & \left(n/(4k^3) - k \right)^{k} \\
&> & \left(n/(8k^3)\right)^{k} \quad \mbox{ (since $n\ge 8k^4$)}\\
&> & (k+1)  \alpha n^{k}  \quad \mbox{ (since $\alpha< ((k+1)8^{k}k^{3k}))^{-1}$}.
\end{eqnarray*}
This implies that there exists $i \in [k+1]$ such that
$|N_{\mathcal{F}_{t+1}(k,n)}(u_i) \setminus N_{H}(u_i)| > \alpha n^{k}$,
contradicting the fact that all $u_i$ are $\alpha$-good with respect to $\mathcal{F}_{t+1}(k,n)$.

Therefore, $H$ has a matching $M$ of size $t$.
Moreover, $M \cup \{e_0\}$ is a matching of size $t+1$ unless $e_0=\emptyset$ in which case $H$ is a subgraph of $\mathcal{F}_{t+1}(k,n;t)$.
 \qed

\medskip

We can now prove Theorem~\ref{main} when $\mathcal{F}^{t+1}(k,n)$ is  $\varepsilon$-close to $\mathcal{F}_{t+1}(k,n;t)$.
\begin{lemma}\label{close-lem}
Let $k\ge 3$, $t\ge 1$ and $n$  be integers, and let $\varepsilon, \zeta$ be real numbers such that
 $0 < \varepsilon \ll \zeta < 1$ and $\varepsilon\ll 1/k$,
$t \ge 48k^2$ and $t + 1< (1-\zeta)(1-k(k+1)\sqrt{\varepsilon}) n/k$.
Let ${\cal F}=\{F_1, \ldots, F_{t+1}\}$ be a family of subsets of ${[n]\choose k}$ such that
$|F_i|> {n\choose k}-{n-t\choose k} - {n-t-k \choose k-1} + 1$ for
$i\in [t+1]$.
 Suppose $\mathcal{F}^{t+1}(k,n)$ is  $\varepsilon$-close to $\mathcal{F}_{t+1}(k,n;t)$.
 Then $\mathcal{F}^{t+1}(k,n)$ has a matching of size $t+1$ or $\mathcal{F}^{t+1}(k,n)$ is a subgraph of  $\mathcal{F}_{t+1}(k,n;t)$.
\end{lemma}

\pf We may assume $n\le  2k^3t$ as otherwise the assertion follows from Corollary~\ref{thm:exact2}.
Let $X, [n]$ be the partition classes of ${\cal F}_{t+1}(k,n;t)$, and
let   $X := \{x_1, x_2, ..., x_{t+1}\}$.
Note that each edge of ${\cal F}_{t+1}(k,n;t)$ intersects $[t]$.


Let $B$ denote the set of vertices in $\mathcal{F}^{t+1}(k,n)$ that
are not $\sqrt{\varepsilon}$-good with respect to  $\mathcal{F}_{t+1}(k,n;t)$. Since $\mathcal{F}^{t+1}(k,n)$ is
$\varepsilon$-close to $\mathcal{F}_{t+1}(k,n;t)$, {$|B|\leq (k+1)(1+1/k)^{k+1}\sqrt{\varepsilon}n\leq 4(k+1)\sqrt{\varepsilon}n$.}
Let $b:=\max\{|B\cap X|, |B\cap [t]|\}$; so $b \le 4(k+1)\sqrt{\varepsilon}n$.
We choose $X_1 \subseteq X, W_1 \subseteq [t]$ such that $B\cap X \subseteq X_1 $, $B\cap [t]\subseteq W_1$, and $|X_1|=|W_1|=b$.

For each $x_{i} \in X \setminus X_1$, $x_i$ is
$\sqrt{\varepsilon}$-good with respect to $\mathcal{F}_{t+1}(k,n;t)$
and we let
$\mathcal{F}_i=\mathcal{F}^{t+1}(k,n)[X_1\cup W_1\cup U \cup \{x_{i}\}]$.
For every $x\in X_1 \cup \{x_{i}\}$, there exists $j\in [t+1]$ such
that $x=x_j$, and we have
\begin{equation*}
\begin{split}
|N_{\mathcal{F}_i}(x)|
&\geq {|N_{\mathcal{F}^{t+1}(k,n)}(x_j)|}-\left({n\choose k}-{n-|[t]\setminus W_1|\choose k}\right) \\
&=|F_j|-\left({n\choose k}-{n-(t-b)\choose k}\right) \\
&>{n-(t-b)\choose k}-{n-t\choose k}-{n-t-k \choose k-1} + 1 \\
&={n-(t-b)\choose k}-{(n-(t-b))-b \choose k} -{(n-(t-b)) - b -k \choose k-1} + 1.
\end{split}
\end{equation*}
Since $n-(t-b) > n/2 \ge 2k^3(k+1)\sqrt{\varepsilon}n > 2k^3 b$,  it follows from
Corollary \ref{thm:exact2} that the family $\{N_{\mathcal{F}_i}(x): x\in
X_1 \cup \{x_{i}\}\}$  admits a rainbow matching of size $b+1$,
or there exists $W^i \in {[n] \setminus ([t]\setminus W_1)\choose b}$
such that $W^i $ is a vertex cover of $N_{\mathcal{F}_i}(x)$ for $x\in
X_1\cup \{x_i\}$.

Suppose there exists $x_i \in X\setminus X_1$ such that $\{N_{\mathcal{F}_i}(x): x\in X_1 \cup \{x_{i}\}\}$  admits a rainbow matching $M_i$ of size $b+1$.
Write $H := \mathcal{F}^{t+1}(k,n)[X \cup [n] \setminus (V(M_i) \cup B)]$.
Let $n' = |V(H) \cap [n]|$.
Since {$b\le 4(k+1)\sqrt{\varepsilon}n$} and $n>48k^2$,  every vertex
in $H$ is $\varepsilon^{1/3}$-good with respect to
$\mathcal{F}_{t-b}(k,n';t-b)$ (with relabeling of vertices so that
$[t]\setminus W_1$ corresponds to $[t-b]$) .
By Lemma \ref{good-lem}, $H$ has a matching $M'$ of size $t-b$. 
Therefore, $M_1 \cup M'$ is a matching in ${\cal F}^{t+1}(k,n)$ of size $t+1$.

Hence, we may assume that for each $x_i\in X\setminus X_1$, there
exists $W^i \in {[n] \setminus ([t]\setminus W_1)\choose b}$ such that $W^i$ is a vertex cover of  $\mathcal{F}_i$.
If $W:=W^i=W^j$ for all $x_i,x_j\in X\setminus X_1$ then $([t] \setminus W_1)\cup W$ is a vertex cover of $F_i$ for all $i\in [t+1]$;  so $\mathcal{F}^{t+1}(k,n)$ is a subgraph of $\mathcal{F}_{t+1}(k,n;t)$.
Hence, we may assume that there exist $i\ne j$ such that  $W^i \setminus W^j \ne \emptyset$.

Now fix $x \in X_1$. Note that each edge in
$N_{\mathcal{F}^{t+1}(k,n)}(x)$ either intersects $([t] \setminus W_1)
\cup (W^i \cap W^j)$ or intersects both $W^i \setminus W^j$ and $W^j
\setminus W^i$.
Let $r:=|W^i \setminus W^j|=|W^j \setminus W^i|$; then
$|([t]\setminus W_1) \cup  (W^i \cap W^j)|= (t-b)+(b-r)=t-r$.  Hence,
$$N_{\mathcal{F}^{t+1}(k,n)}(x) \le \left( {n \choose k} - {n-(t-r) \choose k} \right) + {n-t-r \choose k-2} r^2.$$
Note that
$${n-(t-r) \choose k} - {n-t \choose k} \ge r {n-t \choose k-1} = \frac{(r-1)(n-t)}{k-1} {n-t-1 \choose k-2} + {n-t \choose k-1}.$$
Hence,
\begin{equation*}
\begin{split}
N_{\mathcal{F}^{t+1}(k,n)}(x)
&\le {n \choose k} - {n-t \choose k} - {n-t \choose k-1} + {n-t-1 \choose k-2} r^2 - \frac{(r-1)(n-t)}{k-1} {n-t-1 \choose k-2} \\
&< {n \choose k} - {n-t \choose k} - {n-t-1 \choose k-1}  \text{\quad
  (as $n > 10kb \ge 10kr$ and $t\le n/k$)}\\
&<  {n \choose k} - {n-t \choose k} - {n-t-k \choose k-1}.
\end{split}
\end{equation*}
This leads to a contradiction as $|F_i|> {n\choose k}-{n-t\choose k} - {n-t-k \choose k-1} + 1$ for $i\in [t+1]$. \qed

\section{Non-extremal case}
To deal with the case when  $\mathcal{F}^{t+1}({k},n)$ is not
$\varepsilon$-close to $\mathcal{F}_{t+1}({k},n;t)$, we need two
lemmas proved by the present authors in \cite{LWY20}, both with a slight
variation but same proof.

The first is Lemma 4.2 in \cite{LWY20}.
The only difference between that lemma and the statement below is that  we  replace the assumption
``$d_H(x_i)>{n\choose k}-{n-t+1\choose k}$ for $i\in [t]$''  in Lemma
4.2 of \cite{LWY20} by ``$d_H(x_i)>{n\choose k}-{n-t+1\choose k} -
\rho n^k$ for $i\in [t]$''. We omit the proof as it is almost
identical to the proof of Lemma 4.2 in \cite{LWY20}.

\begin{lemma}\label{Absorb-lem-new}
	Let $n,k,t$ be integers and  $\zeta, \rho$ be real numbers,
        such that $1/n\ll \rho\ll \zeta\ll 1/k$ and $ n/(2k^3) \le t \le (1-\zeta)n/k$.
Let $H$ be a $(1,k)$-partite $(k+1)$-graph with partition classes $\{x_1, \ldots, x_{\lfloor n/k\rfloor} \}, [n]$ such that
$d_H(x_i)>{n\choose k}-{n-t+1\choose k} - \rho n^k$ for $i\in [t]$ and $d_H(x_i)={n\choose k}$ for $i=t+1, \ldots, \lfloor n/k\rfloor$.
Then for any $c$ with $0<c \ll \zeta$,
there exists  a matching $M$ in $H$ 
	such that $|M|\le 2k c n$ and, for any balanced subset $S\subseteq V(H)$ 
       with $|S|\le (k+1)c^{1.5} n/2$,       $H[V(M)\cup S]$ 
         has a perfect matching.
\end{lemma}

The second result we need is Lemma 6.4 in \cite{LWY20}, with the
assumption $n/(3k^2) \le t$ replaced by  $n/(2k^3)\le t$ in the
statement below. Again, the same proof in  \cite{LWY20} works.  Recall
the definition of $\mathcal{H}_t(k,n;t)$.

\begin{lemma}
\label{Span-subgraph-new}
Let $k\ge 3$ be an integer, $0 < c\ll\rho \ll \varepsilon \ll 1$ be real numbers, $n\in k\mathbb{Z}$ sufficiently large, and let $t$ be an integer with $n/(2k^3) \le t \le (\frac{1}{2}+c)n/k$.
Let $H$ be a $(1,k)$-partite $(k+1)$-graph with partition classes $X,[n]$ such that $k|X|=n$.
Let  $A_1$ and $A_2$ be a partition of $A$ such that $|A_1|=t$ and $|A_2|=n/k-t$. Suppose that $d_{H}(x)>{n\choose k}-{n-t+1\choose k}-\rho n^k$ for all $x\in A_1$ and $d_{H}(x)={n\choose k}$ for all $x\in A_2$. If $H$ is not $\varepsilon$-close to $\mathcal{H}_t(k,n;t)$, then there exists a spanning subgraph $H'$ of $H$
such that the following conditions hold:
\begin{itemize}
	 	\setlength{\itemsep}{0pt}
		\setlength{\parsep}{0pt}
		\setlength{\parskip}{0pt}
		\item[$(1)$] For all $x\in V(H')$, with at most $n^{0.99}$ exceptions,
                               $d_{H'}(x)=(1\pm n^{-0.01})n^{0.2}$.
		\item[$(2)$] For all $x\in V(H')$, $d_{H'}(x)< 2 n^{0.2}$.
		\item[$(3)$] For any two distinct $x,y\in V(H')$, $d_{H'}(\{x,y\})< n^{0.19}$.
	\end{itemize}
\end{lemma}

Lemma~\ref{Span-subgraph-new} allows us to apply the following result  attributed to Pippenger \cite{PS}
(see Theorem 4.7.1 in \cite{AS08}). An {\it edge cover} in a hypergraph $H$
is a set of edges whose union is $V(H)$.
\begin{theorem}[Pippenger]
\label{nibble-new}
	For every integer $k\ge 2$ and reals $r\ge 1$ and $a>0$, there are $\gamma=\gamma(k,r,a)>0$ and $d_0=d_0(k,r,a)$ such that for every $n$ and $D\ge d_0$ the following holds: Every $k$-uniform hypergraph $H=(V,E)$ on a set $V$ of $n$ vertices in which all vertices have positive degrees and which satisfies the following conditions:
	\begin{itemize}
	 	\setlength{\itemsep}{0pt}
		\setlength{\parsep}{0pt}
		\setlength{\parskip}{0pt}
		\item[$(1)$] for all vertices $x\in V$ but at most $\gamma n$ of them, $d_H(x)=(1\pm \gamma)D$;
		\item[$(2)$] for all $x\in V$, $d_H(x)<r D$;
		\item[$(3)$] for any two distinct $x,y\in V$, $d_H(\{x,y\})<\gamma D$;
	\end{itemize}
	contains an edge cover of at most $(1+a)(n/k)$ edges.
\end{theorem}


\noindent\textbf{Proof of Theorem \ref{main_uniform}.}
By Corollary \ref{thm:exact2}, we may assume that $2k(t+1) < n< 2k^3 (t+1)$.
Let $0 < \varepsilon  \ll {1/k} < 1$ be sufficiently small and
$t > t_0(k)=\max\{48k^2, n/(2k)^3\}$ be sufficiently large.
By Observation 1, it suffices to show $\mathcal{F}^{t+1}(k,n)$ has a matching of size $t+1$.
Applying Lemma \ref{close-lem} to $\mathcal{F}^{t+1}(k,n)$ (by
choosing $\zeta < 1/3$ for instance),
we may assume that $\mathcal{F}^{t+1}(k,n)$ is not $4k^5 \varepsilon$-close to $\mathcal{F}_{t+1}(k,n;t)$. That is,
$\mathcal{H}^{t+1}(k,n)$ is not $\varepsilon$-close to $\mathcal{H}_{t+1}(k,n;t+1)$ by Observation 2 as $n<2k^3(t+1)$.

Now we apply Lemma \ref{Absorb-lem-new} to $\mathcal{H}^{t+1}(k,n)$
with $t\leq n/(2k)$ and sufficiently small $\zeta$.
Thus there exists some constant $0 < c \ll \varepsilon$ such that
$n-2k^2cn \ge (2k-2k^4c)(t+1)$,  and $\mathcal{H}^{t+1}(k,n)$ contains an absorbing matching $M_1$ with $m_1:=|M_1|\leq 2kc n$ and for any balanced subset $S$ of vertices with $|S|\leq (k+1)c^{1.5} n$, $\mathcal{H}^{t+1}(k,n)[V(M_1)\cup S]$ has a perfect matching. Let $H:=\mathcal{H}^{t+1}(k,n)-V(M_1)$ and $n_1 := n-km_1$.

Next, we see that $H$ is not $(\varepsilon/2)$-close to $\mathcal{H}_{t+1}(k,n-km_1;t+1)$. For, suppose otherwise.
Then
\begin{align*}
& |E(\mathcal{H}_{t+1}(k,n;t+1)) \setminus  E(\mathcal{H}^{t+1}(k,n))| \\
&\le |E(\mathcal{H}_{t+1}(k,n-km_1;t+1)) \setminus E(H)| + |e \in E(\mathcal{H}_{t+1}(k,n;t+1)) : e \cap V(M_1) \neq \emptyset| \\
&\leq (\varepsilon/2) (n+n/k)^{k+1} + 2k(k+1)cn \cdot n^k \\
&\le {\varepsilon (n+n/k)^{k+1}.}
\end{align*}
This is a contradiction as   $\mathcal{H}^{t+1}(k,n)$ is not $\varepsilon$-close to $\mathcal{H}_{t+1}(k,n;t+1)$.

Let {$Q,[n]$} denote the partition classes of $\mathcal{H}^{t+1}(k,n)$,
and let {$X$} consist of all vertices of {$Q$} contained in
$\mathcal{F}^{t+1}(k,n)$.   Let $A_1=X\cap V(H)$ and $A_2=(Q\setminus
X)\cap V(H)$. Then  we have $d_{H}(x)>{n\choose k}-{n-t+1\choose k}-\rho n^k$ for all $x\in A_1$ and $d_{H}(x)={n\choose k}$ for all $x\in A_2$.
Since $n_1 =n-km_1\ge n - 2k^2cn \ge (2k-2 k^4 c)(t+1)$, by Lemma \ref{Span-subgraph-new} $H$ has a spanning subgraph $H_1$ such that
\begin{itemize}
		\item[(1)] for all but at most $n_1^{0.99}$ vertices $x\in V(H_1)$,
                $d_{H_1}(x)=(1\pm n_1^{-0.01})n_1^{0.2}$;
		\item[(2)] for all $x\in V(H_1)$, $d_{H_1}(x)< 2 n_1^{0.2}$;
		\item[(3)] for any two distinct $x,y\in V(H_1)$, $d_{H_1}(\{x,y\})< n_1^{0.19}$.
	\end{itemize}

Hence by applying Lemma \ref{nibble-new} to $H_1$ by choosing $a$ with $0
< a \ll c^{1.5}$,
$H_1$ contains an edge cover of size at most $(1+a)((n_1/k+n_1)/(k+1))$.
Thus, at most $a(n_1/k+n_1)$ vertices are each covered by more than one edge
in the cover. Hence,
after removing at most $a(n_1/k+n_1)$ edges from the edge cover, we obtain a matching $M_2$ covering all but at most $(k+1)a(n_1/k+n_1) \le 3kan_1 \le 3kan$ vertices.

Now we may choose a subset $S$ of $V(H)\setminus V(M_2)$ such that
$|V(H)\setminus (V(M_2)\cup S)|\leq k$ and if $S\ne \emptyset$ then
$|S\cap [n]|={k|S\cap Q|}$.
Since $|S| \le 3kan \le (k+1)c^{1.5}n$, $\mathcal{H}^{t+1}(k,n)[V(M_1)\cup S]$ has a perfect matching, say $M_3$. Thus,
$M_2\cup M_3$ is matching of ${\cal H}^{t+1}(k,n)$ covering all but at most $k$ vertices, and, hence, has size $\lfloor n/k\rfloor$.
Therefore, by Lemma \ref{Rain-PM}, $\mathcal{F}^{t+1}(k,n)$ has a matching of size $t+1$. \qed

\section{Conclusion}

We proved stability results for rainbow matchings for families of
non-uniform hypergraphs. In the uniform case, our result generalizes
the result of Frankl and Kupavskii \cite{FK20} on Conjecture~\ref{FK-stablity}.
We think the following could be true:

\begin{conjecture}
Let  $n,k,t$ be three positive integers such that $n > kt$. Let ${\cal
  F}=\{F_1,\ldots, F_{t+1}\}$ be a  family of subsets of ${[n]\choose
  k}$, such that,  for $1\leq i\leq t+1$,
$$|F_i|> \max\{ e(H_k(n,t,2)), \dots, e(H_k(n,t,k)), e(HM^{(k)}(t,n)) \}.$$
Then ${\cal F}$ admits a rainbow matching of size $t+1$ or there exists $W\in {[n]\choose t}$ such that for all $i\in [t+1]$, $W$ is a vertex cover of $F_i$.
\end{conjecture}



\end{document}